\documentclass[12pt]{amsart}
\usepackage{amsfonts}
\usepackage{amsmath}
\usepackage{amssymb}
\usepackage[margin=1in]{geometry}
\usepackage{tikz}

\newtheorem{theorem}{Theorem}

\newtheorem{proposition}{Proposition}

\theoremstyle{definition}

\newtheorem{remark}{Remark}

\makeatletter
\@namedef{subjclassname@2020}{\textup{2020} Mathematics Subject Classification}
\makeatother

\numberwithin{equation}{section}

\usepackage{mathtools}							
\usepackage{commath}							
\usepackage[overload]{textcase}					
											
\usepackage[english]{babel}						
\usepackage[hidelinks]{hyperref}  					

\newcommand{\N}{\mathbb{N}}
\newcommand{\Z}{\mathbb{Z}}

\newcommand{\R}{\mathbb{R}}
\newcommand{\C}{\mathbb{C}}

\newcommand{\MS}{\mathcal{S}}
\newcommand{\MP}{\mathcal{P}}
\newcommand{\K}{\mathcal{K}}

\newcommand{\A}{\mathcal{A}}
\DeclareMathOperator{\Fp}{F.p.}
\DeclareMathOperator{\supp}{supp}
\DeclareMathOperator{\sgn}{sgn}

\renewcommand{\Re}{\operatorname{Re}}
\renewcommand{\Im}{\operatorname{Im}}

\newcommand{\I}{\mathrm{i}}
\newcommand{\e}{\mathrm{e}}

\newcommand{\vphi}{\varphi}
\newcommand{\eps}{\varepsilon}

\begin{document}

\title[The Fourier-Laplace transforms of certain oscillatory functions]
{An asymptotic analysis of the Fourier-Laplace transforms of certain oscillatory functions}

\author[F.~Broucke]{Frederik Broucke}
\thanks{F. Broucke was supported by the Ghent University BOF-grant 01J04017}

\author[G.~Debruyne]{Gregory Debruyne}
\thanks{G.~Debruyne acknowledges support by Postdoctoral Research Fellowships of the Research Foundation--Flanders (grant number 12X9719N) and the Belgian American Educational Foundation. The latter one allowed him to do part of this research at the University of Illinois at Urbana-Champaign.} 

\author[J.~Vindas]{Jasson Vindas} 
\thanks{J. Vindas was partly supported by Ghent University through the BOF-grant 01J04017 and by the Research Foundation--Flanders through the FWO-grant 1510119N}

\address{Department of Mathematics: Analysis, Logic and Discrete Mathematics\\ Ghent University\\
  Krijgslaan 281\\ B 9000 Ghent\\ Belgium} 
  \email{fabrouck.broucke@UGent.be}
	\email{gregory.debruyne@UGent.be}
  \email{jasson.vindas@UGent.be}


\begin{abstract}
We study the family of Fourier-Laplace transforms
$$
F_{\alpha,\beta}(z)= \Fp  \int_{0}^{\infty} t^{\beta}\exp(\I t^{\alpha}-\I z t)\dif t, \quad \Im z<0,
$$
for $\alpha>1$ and $\beta\in\mathbb{C}$, where Hadamard finite part is used  to regularize the integral when $\Re \beta\leq -1$. We prove that each $F_{\alpha,\beta}$ has analytic continuation to the whole complex plane and determine its asymptotics along any line through the origin. We also apply our ideas to show that some of these functions provide concrete extremal examples for the Wiener-Ikehara theorem and a quantified version of the Ingham-Karamata theorem, supplying new simple and constructive proofs of optimality results for these complex Tauberian theorems.
\end{abstract}

\maketitle

\section{Introduction}

In this article we study the Fourier-Laplace transforms of the family of oscillatory functions
\begin{equation}
\label{Eq: F-L osc: oscillatory functions}
	f_{\alpha, \beta}(t) \coloneqq t^{\beta}\exp(\I t^{\alpha}), \quad t>0,
\end{equation}
where $\alpha>1$ and $\beta\in\C$. When $\Re \beta>-1$, these functions are locally integrable and their Fourier-Laplace transforms are given by
\begin{equation}
\label{Eq: F-L osc: oscillatory functions F-L}
F_{\alpha,\beta}(z)\coloneqq   \int_{0}^{\infty} t^{\beta}\exp(\I t^{\alpha}-\I z t)\dif t \qquad \mbox{for } \Im z<0.
\end{equation}
One can extend the definition of $F_{\alpha,\beta}$ to include any  value $\beta\in\mathbb{C}$ if one considers the Hadamard finite part of the integrals in \eqref{Eq: F-L osc: oscillatory functions F-L} (cf. Subsection \ref{subsec: F-L osc: Preliminaries regularizations}). Note that some instances of the parameters $\alpha$ and $\beta$ lead to well-studied classical functions, such as the Gaussian and the Airy function essentially corresponding to $\beta=0$ and $\alpha=2$ or 3, respectively.

We shall show here that the Fourier-Laplace transforms $F_{\alpha, \beta}$ admit analytic continuation as entire functions and our main aim is to determine their asymptotics throughout the complex plane. Although the study of these questions is of significant intrinsic interest, let us point out that they naturally arise in several applications. 

When $\alpha\geq2$ is an integer and $\beta<-1$, the functions \eqref{Eq: F-L osc: oscillatory functions} and their Fourier transforms, namely, the extensions of \eqref{Eq: F-L osc: oscillatory functions F-L} to the real axis,  naturally occur in the study of multifractal properties of various lacunary Fourier series. Indeed, the first order asymptotics of  \eqref{Eq: F-L osc: oscillatory functions F-L} on the real line in combination with the Poisson summation formula play a crucial role \cite{Chamizo-Ubis2007} in the determination of the pointwise H\"older exponent at the rationals of the family of Fourier series $R_{\alpha,\beta}(x)=\sum_{n=1}^{\infty} n^{\beta}\exp(2\pi \I n^{\alpha} x )$, which are generalizations of  Riemann's classical function \cite[Chapter 7]{ChoimetQueffelec2015}. The results of the present paper might certainly be used to further refine \cite[Theorem 2.1]{Chamizo-Ubis2007}  and exhibit full trigonometric chirp expansions for $R_{\alpha,\beta}$ at each rational point. In \cite{DebruyneSchlage-PuchtaVindas2016}, the asymptotic behavior of $F_{\alpha,-1}$ on the real axis has been determined and has shown to be useful in the construction of concrete instances of Beurling prime number systems for the comparison of abstract prime number theorems. 

As a new application of the functions $f_{\alpha,\beta}$, we now explain how they (and some other close relatives) can be used to establish some optimality results in Tauberian theory. Quantified Tauberian theorems have many important applications in several diverse areas of mathematics, ranging from number theory to operator theory. Accordingly, this kind of theorems has been extensively studied over the past decades. We shall consider a variant of the following model theorem, which is a quantified version of the celebrated Ingham-Karamata Tauberian theorem \cite{ingham1935,karamata1934} (cf. \cite[Chapter 3]{ChoimetQueffelec2015} and \cite[Chapter III]{korevaarbook}).

\begin{theorem}[{\cite[Proposition 3.2]{BorichevTomilov2010}}] \label{th: F-L osc: I-K remainder} 
Let $\tau \in L^{\infty}$ and $\kappa > 0$. Suppose that the Laplace transform $\mathcal{L}\{\tau;s\} = \int^{\infty}_{0} \tau(t) \e^{-st} \mathrm{d}t$ admits an analytic continuation to 
\[
	\Omega = \{s: \Re s \geq - C/ (1+\abs{\Im s})^{\kappa}\},
\]
where it has at most polynomial growth. Then,
\begin{equation}
\label{eq: F-L osc: remainder}
	\int^{x}_{0} \tau(t) \mathrm{d}t = \mathcal{L}\{\tau; 0\} + O\biggl(\left(\frac{\log x}{x}\right)^{1/\kappa}\biggr).
\end{equation}
\end{theorem}
 
In fact, motivated by applications in partial differential equations, the above theorem has recently seen numerous generalizations \cite{BattyBorichevTomilov2016,BattyDuyckaerts2008}. The most general one in terms of the region of analytic continuation and the growth inside such a region is currently given in \cite{Stahn2018}. A natural question is then whether the error term in \eqref{eq: F-L osc: remainder} is sharp. Concerning this question of optimality in Theorem \ref{th: F-L osc: I-K remainder}, two rather different approaches are thus far known. The first one arose in \cite{BorichevTomilov2010} and consisted in a delicate function-theoretic construction, where exactly the optimality of Theorem \ref{th: F-L osc: I-K remainder} was proved. The technique was then refined to show optimality for more general versions of Theorem \ref{th: F-L osc: I-K remainder} in \cite{BattyBorichevTomilov2016} and the most general optimality results achieved via this technique can be found in \cite{Stahn2018}. The second approach only appeared very recently in \cite{DebruyneSeifert1} and crucially depends on a careful application of the open mapping theorem, see also \cite{DebruyneSeifert2} for the most general results obtained by this method. The question then remains whether one can find ``simple" functions showing optimality results. Indeed, the first approach gives a rather non-explicit complicated function, whereas the second functional analysis approach does not even construct a example, it merely shows the existence of one. 

We wish to indicate here that one can indeed find such ``simple" functions, in this way effectively providing a new third approach for addressing optimality questions. Our focus here lies on the simplicity of the functions and any attempt to generality is beyond the scope of this article. Furthermore, we shall not directly answer the optimality question for Theorem \ref{th: F-L osc: I-K remainder}, but for a slightly differently formulated version, although the interested reader may verify via the same techniques developed in this paper that the function $\tau(x) = \exp(\I x^{1+1/\kappa}/ \log^{1/\kappa} x)$ if $x\ge \e$, $\tau(x) = 0$ if $x<\e$, satisfies the hypotheses of Theorem \ref{th: F-L osc: I-K remainder} yet 
\[
	\int_{0}^{x}\tau(t)\dif t = \mathcal{L}\{\tau; 0\} + \frac{\exp\bigl(\I x (x/\log x)^{1/\kappa}\bigr)}{\I(1+1/\kappa)}\left(\frac{\log x}{x}\right)^{1/\kappa} + O\biggl(\frac{\log^{1/\kappa - 1} x}{x^{1/\kappa}}\biggr).
\]

We shall show the optimality of Theorem \ref{th: F-L osc: I-K remainder} where the region of analytic continuation $\Omega$ is altered to 
\[
	\Omega = \left\{s: \Re s \geq - \frac{C \log (\, \abs{\Im s} + 2)} {(1+\abs{\Im s})^{\kappa}}\right\}.
\]
The error term in \eqref{eq: F-L osc: remainder} then becomes $O(x^{-1/\kappa})$. The function $\tau(x) = f_{1+1/\kappa, 0}(x) = \exp(\I x^{1+1/\kappa})$ provides then an extremal example for this theorem. Indeed, $\tau$ is clearly bounded, and the polynomial bounds on the analytic extension of the Laplace transform in $\Omega$ follow from our results in Section \ref{sec: F-L osc: further bounds}. However, by integration by parts one can see that
\[
	\int_{0}^{x}\tau(t)\dif t = \mathcal{L}\{\tau; 0\} + \frac{\exp\bigl(\I x^{1+1/\kappa}\bigr)}{\I(1+1/\kappa)x^{1/\kappa}} + O\bigl(x^{-1-2/\kappa}\bigr).
\]

As a final application to Tauberian theory, we now discuss a recent conjecture by M\"uger on the Wiener-Ikehara Tauberian theorem, another cornerstone in Tauberian theory, see \cite[Chapter III]{korevaarbook} for a historical overview of this theorem and some of its applications. In \cite{Muger2018}, it was conjectured that any non-negative non-decreasing function $S$, whose Mellin transform $ \int^{\infty}_{1}S(x) x^{-s-1} \mathrm{d}x$ after subtraction of a pole term $A/(s-1)$ admits an analytic extension to a half-plane $\Re s > \sigma_{0}$ with $0 < \sigma_{0} < 1$, should obey the asymptotic relation $S(x) = Ax + O(x^{\frac{\sigma_{0} + 2}{3}+ \varepsilon})$, for each $\varepsilon > 0$. This conjecture was already disproved in \cite{DebruyneVindas2018}. In fact, in that paper it was even shown that for each positive function $\delta(x) \rightarrow 0$ there are functions $S$ satisfying all requirements, yet $(S(x) - Ax)/(x\delta(x))$ is unbounded. However, the proof depended on the open mapping theorem, and therefore this functional analysis approach is unable to produce any counterexamples. Here, we can give an explicit function disproving the conjecture of M\"uger, namely, 
$$S(x) = \int^{x}_{1} (1 + \cos(\log^{\alpha} u)) \dif u \quad \mbox{for }  \alpha > 1.$$ Indeed, a simple computation shows that its Mellin transform equals $1/(s-1) - 1/s +(F_{\alpha,0}(i(1-s))+ \overline{F}_{\alpha,0}(i(1-\overline{s})))/(2s)$ which after subtraction of the pole at $1$ admits an analytic extension to the half-plane $\Re s > 0$ due to the fact shown in this article that $F_{\alpha,0}$ is an entire function. On the other hand, integrating by parts yields 
\[
	S(x) = x + \frac{x\sin\bigl(\log^{\alpha} x\bigr)}{\alpha\log^{\alpha-1}x} + O\biggl(\frac{x}{\log^{2(\alpha - 1)}x}\biggr).
\]
We also mention that it is possible to give a constructive proof of the more general negative result from \cite{DebruyneVindas2018} that holds for any function $\delta(x) \rightarrow 0$; see our forthcoming article \cite{Broucke-Debruyne-VindasW-I2020}.

Let us now return to the main subject of our paper, the Fourier-Laplace transforms $F_{\alpha, \beta}$. In Section \ref{sec: F-L osc: analytic continuation} we show they are entire functions; we will actually provide an explicit formula for their analytic continuations.  We point out that the entire extension of $F_{\alpha,\beta}$ for $\beta = -1$ was already communicated in \cite[Theorem 3.1 (a)]{DebruyneSchlage-PuchtaVindas2016}, but the proof given therein was wrong\footnote{The estimates \cite[Eq. (3.12) and (3.13)]{DebruyneSchlage-PuchtaVindas2016} are unclear for $\sigma > 0$. However, upon replacing entire extension by $C^{\infty}$-extension to $\Re s\geq 1$  in \cite[Theorem 3.1(a)]{DebruyneSchlage-PuchtaVindas2016} the statement and proof become correct. Since the entire extension was not used elsewhere in that article, the main results of \cite{DebruyneSchlage-PuchtaVindas2016} are not compromised. }. Section \ref{sec: F-L osc: asymptotics} is devoted to an asymptotic analysis of $F_{\alpha,\beta}$. We shall obtain full asymptotic series on any line through the origin. The asymptotic behavior will display Stokes phenomenon, having qualitatively different asymptotic behavior on the two sectors $\{z:\: -\pi-\pi/\alpha<\arg z<0\}$ and $\{z:\: 0<\arg z<\pi-\pi/\alpha\}$; see the asymptotic formulas \eqref{eq: F-L osc: exp case 1} and \eqref{eq: F-L osc: exp case 2}, respectively. On their boundary rays, the asymptotic behavior will essentially be a mixture of the previous two cases. When $z=x$ is real and positive for example, we have the following asymptotic series.
\begin{proposition}
\label{FL: proposition F-L real}
There are constants $c_{n,\alpha, \beta}$, $d_{n,\alpha,\beta}\in\C$ such that, as $x\to \infty$,
\begin{multline} \label{eq: F-L osc: z real pos.}
	F_{\alpha,\beta}(x) + \sum_{\mathclap{\substack{m,n \\ \beta+n\alpha+m+1=0}}}\frac{\I^{n-m}x^{m}}{n!\,m!}(\log x + \pi\I/2) \sim \sum_{n=0}^{\infty}\frac{c_{n,\alpha,\beta}}{x^{\beta+n\alpha+1}}\\
	{} + \exp\bigl(-\I \alpha^{-1/(\alpha-1)}(1-1/\alpha)x^{\frac{\alpha}{\alpha-1}}\bigr)x^{\frac{\beta+1-\alpha/2}{\alpha-1}}\sum_{n=0}^{\infty}\frac{d_{n,\alpha,\beta}}{x^{n\alpha/(\alpha-1)}}.
\end{multline}
The coefficients $c_{n,\alpha,\beta}$ and $d_{n,\alpha, \beta}$ are given by \eqref{eq: F-L osc: c_{n,alpha,beta}} and \eqref{eq: F-L osc: d_{n,alpha,beta}}, respectively.
\end{proposition}

It should be noted that the main leading terms of the asymptotic expansion \eqref{eq: F-L osc: z real pos.} of Proposition \ref{FL: proposition F-L real} were essentially obtained in \cite{Chamizo-Ubis2007,DebruyneSchlage-PuchtaVindas2016} in some cases by employing Littlewood-Paley decompositions of the unity. 
 Our approach in this article is based on different technology. We exploit here the Estrada-Kanwal moment asymptotic expansion \cite{EstradaKanwal1990,EstradaKanwalbook} in combination with contour integration to deduce asymptotic series expansions. This technique turns out to provide a unified way to deal with the distinct cases of asymptotic behavior that we shall encounter in Section \ref{sec: F-L osc: asymptotics}; in addition, it directly yields desired uniformity of the asymptotic expansions on closed subsectors.  Finally, the article concludes with some polynomial bounds in Section \ref{sec: F-L osc: further bounds} for $F_{\alpha,\beta}$ on hourglass-shaped neighborhoods of the real line.

\section{Preliminaries}
\label{sec: F-L osc: Preliminaries}
We collect in this section several useful notions that play a role in the article. We always write $z=x+\I y$ for complex variables. We will often employ Vinogradov's notation $f(t)\ll g(t)$ to denote the asymptotic bound $f(t)=O(g(t))$. We shall make extensive use of Schwartz  distribution calculus in our manipulations; see the monographs \cite{EstradaKanwalbook,PilipovicStankovicVindas2012} for background material on asymptotic analysis in the distributional setting. Our notation for distributions is as in \cite{EstradaKanwalbook}. In particular, we denote the dual pairing between a distribution $f$ and a test function $\psi$ as $\langle f, \psi \rangle$, or as $\langle f(t), \psi(t) \rangle$ with the use of a dummy  variable of evaluation.

\subsection{Distributional regularization}
\label{subsec: F-L osc: Preliminaries regularizations}
Let us first clarify our interpretation of $f_{\alpha, \beta}(t) = t^{\beta}\exp(\I t^{\alpha})$ as Schwartz distributions. We take them as 0 on $(-\infty,0)$. If $\Re \beta > -1$, then $f_{\alpha, \beta}  \in L^{1}_{\mathrm{loc}}$, and since they are of polynomial growth, they can be viewed as elements of the space of tempered distributions $\MS'$. When $\Re \beta \leq -1$, they do not define distributions automatically and we have to consider regularizations. This can be done in many ways (see e.g.\ \cite[Section 2.4]{EstradaKanwalbook}), but it is desirable that the property $tf_{\alpha,\beta}(t) = f_{\alpha,\beta+1}(t)$ remains true. (Then $F_{\alpha, \beta}' = -\I F_{\alpha, \beta+1}$.) 
This is the case when we regularize them by taking Hadamard finite part \cite[page 67]{EstradaKanwalbook}. Suppose $\psi \in \MS$. Define
\begin{multline*}
	\langle f_{\alpha,\beta}, \psi\rangle 	\coloneqq \Fp \int_{0}^{\infty}t^{\beta}\e^{\I t^{\alpha}}\psi(t)\dif t = \int_{1}^{\infty}t^{\beta}\e^{\I t^{\alpha}}\psi(t)\dif t \\
	+ \int_{0}^{1}t^{\beta}\biggl(\e^{\I t^{\alpha}}\psi(t) - \sum_{\mathclap{\substack{m,n \\ m+n\alpha + \Re\beta\leq -1}}} 
		\frac{\I^{n} \psi^{(m)}(0)}{n! \, m!}t^{n\alpha + m}\biggr)\dif t \quad + \quad \sideset{}{'}\sum_{\mathclap{\substack{m,n \\ m+n\alpha + \Re\beta\leq -1}}}\frac{\I^{n} \psi^{(m)}(0)}{n!\, m!}\frac{1}{\beta+n\alpha+m+1}.
\end{multline*}
Here the notation $\sum'$ means that the possible terms with $\beta+n\alpha+m+1=0$ are excluded. This choice for the regularizations also has the property that for fixed $\alpha$ the map $\beta \mapsto f_{\alpha,\beta}$ is meromorphic; it has poles at $\beta = -n\alpha - m -1$, $n, m \in \N$ with residues
\[
	\sum_{\substack{m,n \\ m+n\alpha+\beta+1=0}}\frac{\I^{n}}{n!}\frac{(-1)^{m}\delta^{(m)}}{m!},
\]
where $\delta^{(m)}$ denotes the $m$-th derivative of the Dirac delta distribution.

\subsection{The moment asymptotic expansion}
To obtain asymptotics for the functions $F_{\alpha, \beta}$, we will systematically employ various techniques based on the so-called moment asymptotic expansion \cite[Chapter 3]{EstradaKanwalbook} (cf. \cite{EstradaKanwal1990}).  Let $\A$ be a test function space of smooth functions where dilation and derivatives act continuously and consider the corresponding distribution space $\A'$ (its topological dual). A distribution $f\in \A'$ is said to satisfy the \emph{moment asymptotic expansion} if
\begin{equation}
\label{eq: F-L osc: moment expansion}
	f(\lambda t) \sim \sum_{n=0}^{\infty}\frac{(-1)^{n}\mu_{n}\delta^{(n)}(t)}{n! \lambda^{n+1}}, \quad \text{as } \lambda \to \infty, 
\end{equation}
where the $\mu_{n}$ are the moments of $f$ given by $\mu_{n} \coloneqq \langle f(t), t^{n}\rangle$. 
The relation \eqref{eq: F-L osc: moment expansion} is to be interpreted in the weak sense, meaning that the asymptotic expansion holds after evaluation at each test function $\psi\in\mathcal{A}$. Such distributions are also called \emph{distributionally small} at $\infty$.

We shall actually make use of a slight generalization of  \eqref{eq: F-L osc: moment expansion}. Let $(\alpha_{n})_{n}$ be a sequence of complex numbers  such that $(\Re \alpha_{n})_{n}$ is increasing and tends to infinity. Suppose a test function $\psi$ on $(0,\infty)$ has the following expansion 
\begin{equation}
\label{eq: F-L osc: test function at 0}
	\psi(t) \sim \sum_{n=0}^{\infty}c_{n}t^{\alpha_{n}}, \quad \mbox{as } t\to0^{+}.
\end{equation}
Then, we consider the generalized moment asymptotic expansion
\begin{equation}
\label{eq: F-L osc: generalized moment expansion}
	\langle f(\lambda t), \psi(t)\rangle \sim \sum_{n=0}^{\infty}\frac{c_{n}\mu_{\alpha_{n}}}{\lambda^{\alpha_{n}+1}}, \quad \text{as } \lambda \to \infty,
\end{equation}
where $\mu_{\alpha_{n}}$ are the generalized moments of $f$ given by $\mu_{\alpha_{n}} \coloneqq \langle f(t), t^{\alpha_{n}}\rangle$. The (generalized) moment asymptotic expansion and its validity have been thoroughly investigated by Estrada and Kanwal and we refer to \cite{EstradaKanwalbook} for a complete account on the subject. See also \cite{NeytVindas2019MAE, EstradaYangMAE} for some recent contributions.  

 In this article we shall exploit the fact that \eqref{eq: F-L osc: generalized moment expansion} holds for all elements of the ensuing two distribution spaces on $(0,\infty)$. 
 \begin{itemize}
 \item The test function space $\MP\{t^{\alpha_{n}}\}$ consists of all functions $\psi \in C^{\infty}(0,\infty)$ having asymptotic expansions \eqref{eq: F-L osc: test function at 0} and such that
\[
	\forall k\in \N, \forall c>0\colon 	\quad	\psi^{(k)}(t) = O_{k,c}(\e^{ct}), \quad					\text{as } t \to \infty.
\]
This space becomes a Fr\'echet space via the ensuing family of norms
\[
\norm{\psi}_{n} 	\coloneqq  \sup_{0<t\le 1}\frac{\abs{\psi(t) - \sum_{j=0}^{n}c_{j}t^{\alpha_{j}}}}{\abs[1]{t^{\alpha_{n+1}}}} + \sup_{\mathclap{\substack{t>1 \\ j\leq n}}}\abs[1]{\psi^{(j)}(t)\e^{-\frac{t}{n}}}, \qquad n=0,1,2,\dots .
\]
All distributions $f\in \MP'\{t^{\alpha_{n}}\}$ satisfy the generalized moment asymptotic expansion (cf. \cite[Sections 2.9, 2.10, 3.4]{EstradaKanwalbook}).
 
\item The space $\K\{t^{\alpha_{n}}\}$ consists of all functions $\psi \in C^{\infty}(0,\infty)$ having asymptotic expansions \eqref{eq: F-L osc: test function at 0} and such that
\[	\exists \,q \in \Z, \forall k \in \N\colon\quad	\psi^{(k)}(t) = O_{k}(t^{q-k}), \quad						\text{as } t \to \infty.
\]
This space is topologized as the inductive limit of the spaces $\K_{q}\{t^{\alpha_{n}}\}$ as $q\to\infty$, where $\K_{q}\{t^{\alpha_{n}}\}$ is the Fr\'{e}chet space of functions $\psi$ which satisfy the above conditions for this fixed $q$, provided with the family of norms
\[
\norm{\psi}_{n, q} \coloneqq 
 \sup_{0<t\le 1}\frac{\abs{\psi(t) - \sum_{j=0}^{n}c_{j}t^{\alpha_{j}}}}{\abs[1]{t^{\alpha_{n+1}}}} +\sup_{\mathclap{\substack{t>1 \\ j\leq n}}}\abs[1]{\psi^{(j)}(t)t^{j-q}}, \qquad n=0,1,2,\dots .
\]
Every distribution $f\in \K'\{t^{\alpha_{n}}\}$ satisfies the generalized moment asymptotic expansion (cf. \cite[Sections 2.9, 2.10, and 3.7]{EstradaKanwalbook}).
	
\end{itemize}

\section{The analytic continuation}
\label{sec: F-L osc: analytic continuation}

We now show that the Fourier-Laplace transform of $f_{\alpha, \beta}$, given by 
\[F_{\alpha, \beta} (z) := \langle f_{\alpha, \beta}(t),\e^{-\I z t}\rangle=\Fp \int_{0}^{\infty}t^{\beta}\exp\left(\I t^{\alpha}-\I z t\right)\dif t \quad \mbox{for }y=\Im z<0,  \]
has holomorphic extension to the whole complex plane. For it, we shift the contour of integration to the ray $\arg\zeta = \pi/(2\alpha)$ where $\I \zeta^{\alpha}$ is real and negative\footnote{We use the principal branch of the logarithm.}. By Cauchy's theorem,
\[
	F_{\alpha, \beta}(z) = \hat{f}_{\alpha, \beta}(z) = \Fp \lim_{\eps\to 0} \int_{\Gamma_{1,\eps}\cup\Gamma_{2,\eps}}\zeta^{\beta}\exp(\I \zeta^{\alpha} - \I z\zeta)\dif \zeta,
\]
where $\Gamma_{1,\eps}$ is the arc of the circle of radius $\eps$ and center at the origin between the points  $\eps$ and $\eps \e^{\I \pi/(2\alpha)}$, and $\Gamma_{2,\eps}$ is the half line $\{ \e^{\I\pi/(2\alpha)}t : t\in [\eps, \infty) \}$. Indeed, defining $\Gamma_{R}$ as the circle arc of radius $R$ from $\e^{\I\pi/(2\alpha)}R$ to $R$, it is easy to see that, for $y<0$, $\int_{\Gamma_{R}}\zeta^{\beta}\exp(\I\zeta^{\alpha}-\I z\zeta)\dif\zeta \to 0$ as $R\to\infty$.
After a small computation, one gets the following expression,
\begin{align}
	F_{\alpha,\beta}(z) 	
	&=\e^{\I(\beta+1)\pi/(2\alpha)}\Fp \int_{0}^{\infty}t^{\beta}\exp(-t^{\alpha}-\I \e^{\I\pi/(2\alpha)}zt)\dif t \quad + \quad \sum_{\mathclap{\substack{m,n \\ m+n\alpha+\beta+1=0}}}\frac{\I^{n}(-\I z)^{m}}{n! \, m!}\frac{\pi}{2\alpha}\I \nonumber\\
	&= \e^{\I(\beta+1)\pi/(2\alpha)}\Biggl[ \int_{1}^{\infty}t^{\beta}\exp\bigl(-t^{\alpha}-\I \e^{\I\pi/(2\alpha)}zt\bigr)\dif t \nonumber\\
					&{} + \int_{0}^{1}t^{\beta}\biggl(\exp\bigl(-t^{\alpha}-\I \e^{\I\pi/(2\alpha)}zt\bigr) - \sum_{\mathclap{\substack{m,n \\ m+n\alpha + \Re\beta\leq -1}}}\frac{(-1)^{n}(-\I z)^{m}\e^{\I m\pi/(2\alpha)}}{n!\,m!}t^{n\alpha+m}\biggr)\dif t \nonumber \\
					&{} +  \sideset{}{'}\sum_{\mathclap{\substack{m,n \\ m+n\alpha + \Re\beta\leq -1}}}\frac{(-1)^{n}(-\I z)^{m}\e^{\I m\pi/(2\alpha)}}{n!\,m!}\frac{1}{\beta+n\alpha+m+1} \Biggr] 
					\quad + \quad \sum_{\mathclap{\substack{m,n \\ m+n\alpha+\beta+1=0}}}\frac{\I^{n}(-\I z)^{m}}{n! \, m!}\frac{\pi}{2\alpha}\I. \label{eq: F-L osc: analytic cont}
\end{align}
The right hand side however is well defined for any $z\in\C$ since $\alpha>1$, so this expression yields the desired entire continuation.


\section{Asymptotic expansion on rays}
\label{sec: F-L osc: asymptotics}

Write $z = R\e^{\I\theta}$. We will derive in this section an asymptotic series expansion for $F_{\alpha, \beta}(R\e^{\I\theta})$ as $R\to\infty$. We distinguish three cases for the angle $\theta$: the sector $\{z:\: -\pi-\pi/\alpha<\arg z<0\}$, the sector $\{z:\: 0<\arg z<\pi-\pi/\alpha\}$, and their boundaries. 

\bigskip

\noindent\textbf{Case 1:} $-\pi-\pi/\alpha<\theta<0.$

In this case we have the following expansion for $F_{\alpha, \beta}(R\e^{\I\theta})$, uniformly on closed subsectors:
\begin{multline}
	\label{eq: F-L osc: exp case 1}
	F_{\alpha, \beta}(R\e^{\I\theta}) + \sum_{\mathclap{\substack{m,n \\ \beta+n\alpha+m+1=0}}}\frac{\I^{n}\e^{\I (\theta-\pi/2)m}R^{m}}{n!\,m!}(\log R + \bigl(\theta+\pi/2\bigr)\I) \\
	\sim \sum_{n=0}^{\infty} \frac{\exp\bigl(\I(n\pi/2 - (\theta+\pi/2)(\beta+n\alpha+1))\bigr)\Gamma^{\ast}(\beta+n\alpha+1)}{n!R^{\beta+n\alpha+1}},
\end{multline}

\smallskip

\noindent where $\Gamma^{\ast}(z)$ equals the Euler gamma function $\Gamma(z)$ when $z\notin -\N$, and otherwise Hadamard finite part values are used in the case that $z \in -\N$, that is,
\[
	\Gamma^{\ast}(-k) \coloneqq \Fp \int_{0}^{\infty}\e^{-t}t^{-k-1}\dif t = \frac{(-1)^{k}}{k!}\bigg(-\gamma+ \sum _{j=1}^{k}\frac{1}{j}\bigg).
\]
Here $\gamma$ is the Euler-Mascheroni constant.

In order to deduce \eqref{eq: F-L osc: exp case 1}, we consider two overlapping subcases: the case where $-\pi<\theta<0$ (equivalently, $y<0$), and the case where $-\pi-\pi/\alpha<\theta<-\pi/\alpha$. 

In the first subcase, we have the following expression for $F_{\alpha,\beta}$ (which is the original form of the Fourier transform, before shifting the contour):
\begin{equation}
\label{eq: F-L osc: F for y<0}
	F_{\alpha, \beta}(R\e^{\I\theta}) = \Fp \int_{0}^{\infty}t^{\beta}\exp\bigl(\I t^{\alpha} + R\e^{\I(\theta-\pi/2)}t\bigr)\dif t.
\end{equation}
The idea is now to relate the above expression to an evaluation $\langle g(Rt), f_{\alpha, \beta}(t) \rangle$ for a distribution $g$ which is distributionally small at infinity, so that $g$ satisfies the moment asymptotic expansion. Making this precise, we consider the space $\MP\{t^{\beta+n\alpha}\}$ from Section \ref{sec: F-L osc: Preliminaries} (with $\alpha_{n} = \beta + n\alpha$). Clearly, we have that $f_{\alpha,\beta} \in \MP\{t^{\beta+n\alpha}\}$, with $c_{n} = \I^{n}/n!$ in its asymptotic expansion \eqref{eq: F-L osc: test function at 0} as $t\to0^{+}$. The distribution $g=g_{\theta} \in \MP'\{t^{\beta+n\alpha}\}$ will be defined as a regularization of the function $\exp(\e^{\I(\theta-\pi/2)}t)$. If $-\pi<\theta<0$, then $\cos(\theta-\pi/2)<0$ so that $\exp(\e^{\I(\theta-\pi/2)}t)\psi(t)$ is integrable away from the origin for every test function $\psi\in\MP\{t^{\beta+n\alpha}\}$; this product might however not be integrable near the origin. We choose the regularization corresponding to the expression \eqref{eq: F-L osc: F for y<0}. For $\psi \in  \MP\{t^{\beta+n\alpha}\}$,
\begin{multline*}
	\langle g(t), \psi(t) \rangle	\coloneqq \Fp \int_{0}^{\infty}\exp\bigl(\e^{\I(\theta-\pi/2)}t\bigr)\psi(t)\dif t = \int_{1}^{\infty}\exp\bigl(\e^{\I(\theta-\pi/2)}t\bigr)\psi(t)\dif t \\
						 +  \int_{0}^{1}\biggl(\exp\bigl(\e^{\I(\theta-\pi/2)}t\bigr)\psi(t) - \sum_{\mathclap{\substack{m,n \\ m+n\alpha + \Re\beta\leq -1}}} c_{n}\frac{\e^{\I(\theta-\pi/2)m}}{m!}t^{\beta+ n\alpha+m}\biggr)\dif t 
						 + \sideset{}{'}\sum_{\mathclap{\substack{m,n \\ m+n\alpha + \Re\beta\leq -1}}}c_{n}\frac{\e^{\I(\theta-\pi/2)m}}{m!}\frac{1}{\beta+n\alpha+m+1}.
\end{multline*}
This defines a continuous linear functional on $\MP\{t^{\beta+n\alpha}\}$, and one readily sees that 
\[
	\langle g(Rt), f_{\alpha, \beta}(t) \rangle = F_{\alpha,\beta}(R\e^{\I\theta}) + 
	\sum_{\mathclap{\substack{m,n \\ m+n\alpha+\beta+1=0}}}\frac{\I^{n}\e^{\I(\theta-\pi/2)m}R^{m}}{n!\,m!}\log R.
\]
We remark that for fixed $\alpha$ the last sum is non-empty only for countably many values of $\beta$, namely, the poles of the vector-valued meromorphic function $\beta \mapsto f_{\alpha,\beta}$.

One verifies via contour integration that the generalized moments of $g$ are given by
\[
\langle g(t), t^{\beta+n\alpha}\rangle = \e^{-\I(\theta+\pi/2)(\beta+n\alpha+1)}\Gamma^{\ast}(\beta+n\alpha+1) - \delta_{m, -n\alpha-\beta-1}\frac{\e^{\I(\theta-\pi/2)m}}{m!}(\theta+\pi/2)\I.
\]
Here $\delta_{m, -n\alpha-\beta-1}$ stands for the Kronecker delta, that is, 1 if $-n\alpha-\beta-1$ equals the nonnegative integer $m$, and 0 otherwise.
Since $g$ satisfies the generalized moment asymptotic expansion \eqref{eq: F-L osc: generalized moment expansion}, we readily obtain the expansion \eqref{eq: F-L osc: exp case 1}.
Upon inspecting the error terms in such an expansion\footnote{See e.g. \cite[Eq. (3.41), p.~116]{EstradaKanwalbook} for an explicit expression of the error term, which carries over to other distribution spaces where the generalized moment asymptotic expansion holds (cf. \cite[Sections 3.4 and 3.7]{EstradaKanwalbook}). The error terms only depend on a dual seminorm of the $g_{\theta}$ and they are uniformly bounded in the ranges under consideration.}, one sees that they are uniform when $-\pi+\eps\le \theta \le - \eps$, with arbitrary $\eps>0$.

The second subcase is similar, but we start from a different expression for $F_{\alpha,\beta}$. In \eqref{eq: F-L osc: analytic cont} we rotate the contour of integration once again over an angle $\pi/(2\alpha)$, and after some computations one gets the following expression for $F_{\alpha, \beta}$,
\[
	F_{\alpha, \beta}(R\e^{\I\theta}) = \e^{\I\pi(\beta+1)/\alpha}\Fp \int_{0}^{\infty}t^{\beta}\exp\bigl(-\I t^{\alpha} + \e^{\I(\theta-\pi/2+\pi/\alpha)}Rt\bigr)\dif t
	\quad + \quad \sum_{\mathclap{\substack{m,n \\ m+n\alpha+\beta+1=0}}}\frac{\I^{n}R^{m}\e^{\I (\theta-\pi/2)m}}{n!\, m!}\frac{\pi}{\alpha}\I.	
\]
One can now proceed in the same way as in the discussion of the first subcase, and one again finds the expansion \eqref{eq: F-L osc: exp case 1}. So, we have established that this asymptotic series expansion holds in the range $-\pi-\pi/\alpha < \theta <0$ with uniformity on closed subsectors.

\bigskip

\noindent\textbf{Case 2:} $0<\theta<\pi-\pi/\alpha$. 

In this case we get the asymptotic series \eqref{eq: F-L osc: exp case 2} stated below. The first order approximation is 
\begin{equation}
\label{eq: F-L osc: exp case 2 first order}
	F_{\alpha, \beta}(R\e^{\I\theta}) \sim \e^{\I\eta_{1,\theta}}\sqrt{\frac{2\pi}{(\alpha-1)}}\: \alpha^{\frac{-1/2-\beta}{\alpha-1}} R^{\frac{\beta+1-\alpha/2}{\alpha-1}}\exp\bigl(\e^{\I\eta_{2,\theta}}\alpha^{-1/(\alpha-1)}(1-1/\alpha)R^{\frac{\alpha}{\alpha-1}}\bigr),
\end{equation}
where
\[
	\eta_{1,\theta} \coloneqq \frac{\pi}{4} - \frac{\alpha-2\beta-2}{2(\alpha-1)}\theta, \quad \eta_{2,\theta}\coloneqq \frac{\alpha}{\alpha-1}\theta - \frac{\pi}{2}.
\]
Notice that $\cos(\eta_{2,\theta})>0$ in this case.

We shall use the method of steepest descent (also called saddle-point method) to study this case. This is a classical method to obtain asymptotic expansions of integrals of the form $\int_{\Gamma} \e^{Rf(\zeta)}g(\zeta) \dif \zeta$ as $R\to \infty$, where $\Gamma$ is a contour in some region $\Omega$ where $f$ and $g$ are holomorphic. The basic idea is to shift the contour to one that passes through a saddle point $\zeta_{0}$ of $f$, that is, a point for which $f'(\zeta_{0})=0$. If the new contour passes through this point in such a way that $\Re f$ reaches a maximum at $\zeta_{0}$ on this new contour, then one can use some form of the Laplace asymptotic formula to obtain the asymptotics for the integral. We refer to \cite[Section 3.6]{EstradaKanwalbook} for more details.

Starting from expression \eqref{eq: F-L osc: analytic cont}, we will use the method of steepest descent on the integral from $1$ to $\infty$ and this will give the main contribution; the other two terms are $O(\e^{R}R^{\,\abs{\beta}})$ and $O(R^{\,\abs{\beta}-1})$ respectively, and as we will see they are negligible with respect to the main contribution. Set $\kappa \coloneqq1/(\alpha-1)$ , $\vphi\coloneqq \theta-\pi/2 + \pi/(2\alpha)$ and perform the substitution $t=R^{\kappa}s$ to get 
\[
	 \int_{1}^{\infty}t^{\beta}\exp\bigl(-t^{\alpha} + R\e^{\I\vphi}t\bigr)\dif t 
	= R^{\kappa(\beta+1)}\int_{1/R^{\kappa}}^{\infty}s^{\beta}\exp\bigl( R^{\kappa+1}(\e^{\I\vphi}s-s^{\alpha})\bigr)\dif s.
\]
The function 
\begin{equation}
\label{eq: F-L osc: h case 2}
h(\zeta) \coloneqq \e^{\I\vphi}\zeta - \zeta^{\alpha}
\end{equation}
 is holomorphic in $\C \setminus (-\infty,0]$ and has a saddle point at $\zeta_{0}= \alpha^{-\kappa}\e^{\I\kappa\vphi}$. We shift the contour of integration to $\Gamma = \bigcup_{j}\Gamma_{j}$, where
\begin{align*}
	\Gamma_{1} 	&\coloneqq [R^{-\kappa}, r], \quad \text{some small } r>0; \\
	\Gamma_{2}	&\coloneqq \{r \e^{\I\eta} : \eta \text{ ranging from } 0 \text{ to } \kappa\vphi\}; \\
	\Gamma_{3}	&\coloneqq [r\e^{\I\kappa\vphi}, \rho \e^{\I\kappa\vphi}], \quad \text{some large } \rho; \\
	\Gamma_{4}	&\coloneqq \{\rho \e^{\I\eta} : \eta \text{ ranging from } \kappa\vphi \text{ to } 0\}; \\
	\Gamma_{5}	&\coloneqq [\rho, \infty). 
\end{align*}

The main contribution will come from the integral over $\Gamma_{3}$; for the other integrals we have:
\begin{align*}
	R^{\kappa(\beta+1)}\Bigg(\int_{\Gamma_{1}} + \int_{\Gamma_{2}}\Bigg) 	&\ll \e^{\eps R^{\kappa+1}}; \\
	R^{\kappa(\beta+1)}\Bigg(\int_{\Gamma_{4}} + \int_{\Gamma_{5}}\Bigg)	&\ll \e^{-CR^{\kappa+1}}.
\end{align*}
Here, $\eps$ is a number depending on $r$ which can be made arbitrarily small by choosing $r$ arbitrarily small, and $C$ is a number depending on $\rho$ which can be made positive by choosing $\rho$ sufficiently large. On $\Gamma_{3}$, $\Re h(\zeta)$ reaches its maximum at the saddle point; applying \cite[Eq. (3.172), p.~137]{EstradaKanwalbook} gives\footnote{By convention, a change of variables in the space of analytic functionals is done without taking absolute value of the Jacobian, e.g., 
\[ (-1)^{n}\langle \delta^{(n)}(\psi(z)),f(z) \rangle = \od[n]{}{\omega}\biggl(\frac{f(\psi^{-1}(\omega))}{\psi'(\psi^{-1}(\omega))}\biggr)\bigg\rvert_{\omega=0}.\]}
\begin{gather*}
	R^{\kappa(\beta+1)}\int_{\Gamma_{3}}\zeta^{\beta}\exp\bigl(R^{\kappa+1}(\e^{\I\vphi}\zeta-\zeta^{\alpha})\bigr)\dif \zeta \\
	\sim \I R^{\kappa(\beta+1)}\exp\bigl(\alpha^{-\kappa}(1-1/\alpha)\e^{\I\alpha\kappa\vphi}R^{\kappa+1}\bigr) 
	\sum_{n=0}^{\infty}\frac{(-1)^{n}\Gamma(n+1/2)\langle \delta^{(2n)}\bigl(\sqrt{h(\zeta)-h(\zeta_{0})}\bigr), \zeta^{\beta}\rangle}{(2n)!R^{(\kappa+1)(n+1/2)}}.
\end{gather*}
 The branch of $\sqrt{h(\zeta)-h(\zeta_{0})}$ is chosen here in such a way that $\Im \sqrt{h(\zeta)-h(\zeta_{0})}$ is increasing in a neighborhood of $\zeta_{0}$ on $\Gamma_{3}$. 

Since all the other contributions are of lower order than every term in the above asymptotic series, we have the same asymptotic relation (up to a multiplicative constant) for $F_{\alpha,\beta}$:
\begin{multline}
	\label{eq: F-L osc: exp case 2}
	F_{\alpha, \beta}(R\e^{\I\theta}) \sim \e^{\I((\beta+1)\pi/(2\alpha)+\pi/2)} R^{\frac{\beta+1-\alpha/2}{\alpha-1}}\exp\bigl(\e^{\I\eta_{2,\theta}}\alpha^{-1/(\alpha-1)}(1-1/\alpha)R^{\frac{\alpha}{\alpha-1}}\bigr) \\
		{} \times \sum_{n=0}^{\infty}\frac{(-1)^{n}\Gamma(n+1/2)\bigl\langle \delta^{(2n)}\bigl(\sqrt{h(\zeta)-h(\zeta_{0})}\bigr), \zeta^{\beta}\bigr\rangle}{(2n)!R^{\frac{n\alpha}{\alpha-1}}},
\end{multline}
where $h$ is given by \eqref{eq: F-L osc: h case 2} and $\zeta_{0}= \alpha^{-\kappa}\e^{\I\kappa\vphi}$. The asymptotic expansion \eqref{eq: F-L osc: exp case 2} holds uniformly on closed subsectors.

\bigskip

\noindent\textbf{Case 3:} $\theta = 0$ or $\theta = -\pi-\pi/\alpha$.

When $z$ crosses the rays $\theta=0$ and $\theta=-\pi-\pi/\alpha$, the asymptotic behavior of $F_{\alpha, \beta}(z)$ changes qualitatively from \eqref{eq: F-L osc: exp case 1} to \eqref{eq: F-L osc: exp case 2}. On these rays, the asymptotic behavior will be a combination of both \eqref{eq: F-L osc: exp case 1} and \eqref{eq: F-L osc: exp case 2}. To fix ideas, assume $\theta=0$, $z=R$. The other case is actually treated similarly, as we explain below. We start from the following expression for $F_{\alpha, \beta}$, which can be derived as in Section \ref{sec: F-L osc: analytic continuation}, but now only rotating the contour in the integral from $R^{\kappa}$ to $\infty$ (recall that $\kappa=1/(\alpha-1)$). We have
\begin{align}
	F_{\alpha, \beta}(R)	&= \e^{\I\pi(\beta+1)/(2\alpha)}\int_{R^{\kappa}}^{\infty}t^{\beta}\exp\bigl(-t^{\alpha}-\I \e^{\I\pi/(2\alpha)}Rt\bigr)\dif t  \nonumber \\
					&\quad {} +  \I R^{\kappa(\beta+1)}\int_{0}^{\frac{\pi}{2\alpha}}\e^{\I\eta(\beta+1)}\exp\bigl(R^{\kappa+1}(\I \e^{\I\alpha\eta}  - \I \e^{\I\eta})\bigr) \dif \eta \nonumber \\
					&\quad {} + \int_{1}^{R^{\kappa}}t^{\beta}\exp\bigl(\I t^{\alpha} - \I Rt\bigr)\dif t 
					 + \int_{0}^{1}(\dotso - \dotso) + \sideset{}{'}\sum\dotso \nonumber \\
					 &\eqqcolon I_{1} + I_{2} + I_{3} + I_{4} + S.
\end{align}

We will split the integral $I_{3}$ into four pieces using partitions of the unity. The splitting will be done in two steps. In the first step, we split $I_{3}$ into two pieces $I_{3,a}+I_{3,e}$: consider two functions such that $\phi_{a}+\phi_{e} = 1$ on $[1,R^{\kappa}]$ and $0<\eps<1$ with
\begin{alignat*}{3}
	&\phi_{a}\in C^{\infty}[0,\infty), \quad 				&&\supp \phi_{a} \subseteq [0,1+\eps], \quad 				&&\phi_{a} = 1 \text{ on } [0,1+\eps/2]; \\
	&\phi_{e}\in C^{\infty}(-\infty,R^{\kappa}],\quad		&&\supp \phi_{e} \subseteq [1+\eps/2, R^{\kappa}],\quad		&&\phi_{e} = 1 \text{ on } [1+\eps, R^{\kappa}].
\end{alignat*}

The sum 
$$I_{3,a} + I_{4} + S= \Fp \int_{0}^{\infty}t^{\beta}\exp(\I t^{\alpha})\phi_{a}(t)\exp(-\I Rt) \dif t$$ can be treated analogously as in \textbf{Case 1}, with one modification. It is no longer the case that the distribution $g_{\theta} = \exp(\e^{\I(\theta-\pi/2)}t) = \exp(-\I t)$ belongs to $\MP'\{t^{\beta+n\alpha}\}$. To remedy this, we consider the space $\K\{t^{\beta+n\alpha}\}$ from Section \ref{sec: F-L osc: Preliminaries} (with $\alpha_{n} = \beta + n\alpha$). 
We have that our test function $t^{\beta}\exp(\I t^{\alpha})\phi_{a}(t)$ is indeed an element of $\K\{t^{\beta+n\alpha}\}$, as it has compact support. The function $\e^{-\I t}$ can be regularized to yield an element of $\K'\{t^{\beta+n\alpha}\}$: the divergence at the origin is resolved in the same way as in \textbf{Case 1}, while the divergence of the integral away from the origin is resolved by formally integrating by parts enough times so that one gets an absolutely convergent integral. More precisely, given $\psi\in \K_{q}\{t^{\beta+n\alpha}\}$  with $\psi(t) \sim c_{0}t^{\beta} + c_{1}t^{\beta+\alpha}+\dotsb$, one regularizes the divergent integral $\int_{0}^{\infty} \e^{-\I t}\psi(t)\dif t$ as 
\begin{gather*}
	\int_{0}^{1}\biggl(\e^{-\I t}\psi(t) - \sum_{\mathclap{\substack{m,n \\ m+n\alpha + \Re\beta\leq -1}}} c_{n}\frac{(-\I)^{m}}{m!}t^{\beta+ n\alpha+m}\biggr)\dif t 
	 + \sideset{}{'}\sum_{\mathclap{\substack{m,n \\ m+n\alpha + \Re\beta\leq -1}}}c_{n}\frac{(-\I)^{m}}{m!}\frac{1}{\beta+n\alpha+m+1}\\
	 {} + \sum_{j=0}^{q+1}(-1)^{j+1}\frac{\e^{-\I}}{(-\I)^{j+1}}\psi^{(j)}(1) + (-1)^{q+2}\int_{1}^{\infty}\frac{\e^{-\I t}}{(-\I)^{q+2}}\psi^{(q+2)}(t)\dif t.
\end{gather*}
Using the moment asymptotic expansion on this regularization will give that the asymptotics of $I_{3,a}+I_{4}+S$ are exactly like \eqref{eq: F-L osc: exp case 1} in \textbf{Case 1} with $0$ substituted for $\theta$.

Our second step is to deal with the integral $I_{3,e}$. We first perform the substitution $t=R^{\kappa}s$ to get
\[
I_{3,e}= R^{(\beta+1)\kappa}\int_{R^{-\kappa}}^{1}s^{\beta}\exp\bigl(-\I  R^{\kappa+1} h(s)\bigr)\phi_{e}(R^{\kappa}s)\dif s,
\]
where 
\begin{equation}
\label{eq: F-L osc: h case 3}
h(s) \coloneqq s-s^{\alpha}.
\end{equation}
We will estimate $I_{3,e}$ using the stationary phase principle. The function $h$ has a unique stationary point $s_{0}=\alpha^{-\kappa}$; $h'(s_{0})=0$.  This stationary point is contained in $[R^{-\kappa}(1+\eps), 1]$ provided that $R$ is sufficiently large, say $R>2^{1/\kappa}\alpha$. In order to single out the contributions from the endpoints and the interior stationary point, we further split the integral $I_{3,e}$ into three pieces using $\phi$ with $\phi_{b} + \phi_{c} + \phi_{d} = 1$ on $[0,1]$ and $\eps'$ with $0<\eps'<s_{0}/2$ with
\begin{alignat*}{3}
	&\phi_{b}\in C^{\infty}(\R), \quad		&&\supp \phi_{b} \subseteq (-\infty,s_{0}/2], \quad 		&&\phi_{b} = 1 \text{ on } [0,s_{0}/2-\eps']; \\
	&\phi_{c}\in C^{\infty}(\R), \quad 		&&\supp \phi_{c} \subseteq [s_{0}/2-\eps', 1-\eps'/2], \quad&&\phi_{c} = 1 \text{ on } [s_{0}/2, 1-\eps']; \\
	&\phi_{d}\in C^{\infty}(-\infty,1],\quad			&&\supp \phi_{d} \subseteq [1-\eps', 1],\quad			&&\phi_{d} = 1 \text{ on } [1-\eps'/2, 1].
\end{alignat*}
This yields three integrals $I_{3,e} = I_{3,b}+I_{3,c}+I_{3,d}$; the stationary point $s_{0}$ is contained in the support of $\phi_{c}$ if $\eps'$ is sufficiently small (say $\eps' < (1-s_{0})$). Furthermore, the function $\phi_{e}(R^{\kappa}s)$ is $1$ on the integration intervals of $I_{3,c}$ and $I_{3,d}$ if $R$ is sufficiently large (say $R^{\kappa} \ge (1+\eps)/(s_{0}/2 - \eps')$).

For the integral $I_{3,b}$ we have:
\[
	I_{3,b}=	
			R^{\kappa}\int_{-\infty}^{\infty}(R^{\kappa}s)^{\beta}\exp\bigl(\I R^{\kappa+1}(s^{\alpha}-s)\bigr)\phi_{e}(R^{\kappa}s)\phi_{b}(s)\dif s.
\]
We now show that $I_{3,b} \ll_{n} R^{-n}$ for any $n\in \N$. Perform the substitution $u = h(s)$ and integrate by parts $n$ times to obtain
\[
	I_{3,b} = \frac{R^{\kappa}}{(\I R^{\kappa+1})^{n}}\int_{J}\exp(-\I R^{\kappa+1}u) \od[n]{}{u}\biggl((h^{-1}(u)R^{\kappa})^{\beta}\phi_{b}(h^{-1}(u))\phi_{e}(R^{\kappa}h^{-1}(u))\frac{1}{h'(h^{-1}(u))}\biggr)\dif u,
\]	
where the integration interval is $J=[h(R^{-\kappa}(1+\eps/2)), h(s_{0}/2)]$. On this interval, we have 
\[ 
	\od[j]{}{u}\frac{1}{h'(h^{-1}(u))} \ll_{j} 1+ R^{-\kappa(\alpha-(j+1))} \ll_{j} R^{j\kappa},
\]
so $I_{3,b} \ll_{n} R^{\kappa(1+\,\abs{\beta})}R^{n\kappa}R^{-n(\kappa+1)} = R^{\kappa(1+\,\abs{\beta})-n}$.

The integral $I_{3,c}$ equals
\[
	R^{\kappa(\beta+1)}\int_{s_{0}/2-\eps'}^{1-\eps'/2}\exp\bigl(-\I R^{\kappa+1}h(s)\bigr)s^{\beta}\phi_{c}(s)\dif s.
\]
The integrand is a smooth function whose support is compact and contains the stationary point $s_{0}$. An asymptotic formula can thus be obtained via the stationary phase principle. 
Employing \cite[Eq.~(3.212), p.~146]{EstradaKanwalbook} \footnote{There are some typos there, one should replace $n$ by $2n$ in the phase of the complex exponential and in the factorial, and $n+1$ by $n+1/2$ in the exponent of $\lambda$.} we get
\begin{multline*}
	I_{3,c} \sim R^{\kappa(\beta+1)}\exp(-\I R^{\kappa+1}h(s_{0}))\sum_{n=0}^{\infty}\frac{\exp(\I\pi(2n+1)/4)\Gamma(n+1/2)}{(2n)!R^{(\kappa+1)(n+1/2)}} \\
	{} \times \bigl\langle \delta^{(2n)}\bigl(\sgn(s-s_{0})\sqrt{h(s_{0}) - h(s)}\bigr), s^{\beta}\bigr\rangle.
\end{multline*}

Finally, $I_{3,d}$ will give a contribution from its endpoint $1$, but this will be cancelled by the contribution from the endpoint $0$ of $I_{2}$: $I_{2}+I_{3,d} \ll_{n} R^{-n}$ for every $n\in\N$. Also $I_{1} \ll_{n} R^{-n}$ for every $n\in\N$.

Collecting all terms, we get the asymptotic expansion \eqref{eq: F-L osc: z real pos.}, with 
\begin{align}\label{eq: F-L osc: c_{n,alpha,beta}}
	c_{n,\alpha,\beta} &= \frac{1}{n!}\exp\biggl(-\frac{\I\pi}{2}\bigl(\beta+1+n(\alpha-1)\bigr)\biggr) \Gamma^{\ast}(\beta+n\alpha+1), \\
\label{eq: F-L osc: d_{n,alpha,beta}}
	d_{n,\alpha,\beta} &= \frac{1}{(2n)!}\exp\bigl(\I\pi(2n+1)/4\bigr)\Gamma(n+1/2)\bigl\langle \delta^{(2n)}\bigl(\sgn(s-s_{0})\sqrt{h(s_{0}) - h(s)}\bigr), s^{\beta}\bigr\rangle.
\end{align}
where $h$ is given by \eqref{eq: F-L osc: h case 3} and $s_{0}=\alpha^{-\kappa}$. Explicitly, we have the following expression for $d_{0,\alpha,\beta}$:
\[
	d_{0,\alpha,\beta} = \e^{\I\pi/4}\sqrt{\frac{2\pi}{\alpha-1}}\alpha^{\frac{-1/2-\beta}{\alpha-1}}.
\]

The case $\theta=-\pi-\pi/\alpha$ is similar, but starting from equation \eqref{eq: F-L osc: analytic cont} we rotate the contour from 0 to $R^{\kappa}$ over an additional angle of $\pi/(2\alpha)$, as in the second subcase of \textbf{Case 1}. One gets:

\smallskip

\begin{align*}
	F_{\alpha,\beta}(R\e^{-\I(\pi+\pi/\alpha)}) 	&+ \sum_{\mathclap{\substack{m,n \\ \beta+n\alpha+m+1=0}}}\frac{\I^{n}\exp(-\I m(3\pi/2 + \pi/\alpha))}{n!\, m!}R^{m}(\log R - \I(\pi/2+\pi/\alpha)) \\
								&\sim \sum_{n=0}^{\infty}\frac{\exp\bigl(\I(n\pi/2 + (\pi/2+\pi/\alpha)(\beta+n\alpha+1)\bigr)\Gamma^{\ast}(\beta+n\alpha+1)}{n! R^{\beta+n\alpha+1}} \\
								&{} + \e^{\I\pi(\beta+1)/\alpha}\exp\bigl(\I \alpha^{-1/(\alpha-1)}(1-1/\alpha)R^{\frac{\alpha}{\alpha-1}}\bigr)R^{\frac{\beta+1-\alpha/2}{\alpha-1}} \\
								& \times  \sum_{n=0}^{\infty}\frac{\exp(-\I\pi(2n+1)/4)\Gamma(n+1/2)}{(2n)!R^{\frac{n\alpha}{\alpha-1}}}\bigl\langle \delta^{(2n)}\bigl(\sgn(s-s_{0})\sqrt{h(s_{0}) - h(s)}\bigr), s^{\beta}\bigr\rangle.
\end{align*}

\section{Bounds in an hourglass-shaped region near the real line} 
\label{sec: F-L osc: further bounds}

In this last section we deduce polynomial bounds for $F_{\alpha, \beta}$ in an hourglass-shaped region near the real axis. Given $C>0$, consider the closed region
\[
	\Omega_{C} \coloneqq \biggl\{ z=x+\I y \in \C: \:  |y| \le C\frac{\log(2+\abs{x})}{(1+\abs{x})^{\kappa}}\biggr\},
\]
where $\kappa=1/(\alpha-1)$. 

For $x$ negative and sufficiently large in absolute value, $z=x+\I y$ lies in the sector treated in \textbf{Case 1}, and by the uniformity of the expansions \eqref{eq: F-L osc: exp case 1} there, we have
\[
	F_{\alpha, \beta}(z) \ll 
	\begin{dcases*}
		\abs{x}^{-1-\Re\beta}			&if $-1-\beta\notin\N$,\\
		\abs{x}^{-1-\beta}\log\abs{x}	&if $-1-\beta\in\N$, 
	\end{dcases*}		
		\qquad \text{when } z\in \Omega_{C} \mbox{ and } x\le0.
\]

When $x$ is positive, we use a similar contour as in \textbf{Case 3}: set $\rho \coloneqq Ax^{\kappa}$ for a parameter $A$ (to be determined below) and rotate the contour in the integral from $\rho$ to $\infty$. We keep $x>1$. For the ``rotated" integral we have
\[
	\int_{\rho}^{\infty}t^{\beta}\exp\bigl(-t^{\alpha} + \e^{-\I\frac{\pi}{2}\left(1-\frac{1}{\alpha}\right)}(x+\I y)t\bigr)\dif t \ll \int_{\rho}^{\infty}t^{\beta}\exp\bigl(-t(t^{\alpha-1} - x - |y|)\bigr)\dif t \ll \e^{-\rho},
\]
since $t^{\alpha-1}-x-|y| \ge A^{\alpha-1}x -x -C(x+1)^{-\kappa}\log( x+2) \ge 2$ if $A>1$ and $x$ is sufficiently large. For the integral over the circle arc we have (using the bounds $2\eta/\pi\le\sin\eta\le\eta$ for $0\le\eta\le\pi/2$)
\begin{multline*}
	\rho^{\beta+1}\int_{0}^{\frac{\pi}{2\alpha}}\e^{\I\beta\eta}\exp\bigl(\I\rho^{\alpha}\e^{\I\alpha\eta} - \I(x+\I y)\rho \e^{\I\eta}\bigr)\I \e^{\I\eta}\dif\eta \\
	\begin{aligned}
		&\ll_{A} x^{\kappa(\Re\beta+1)}\exp(\rho |y|)\int_{0}^{\frac{\pi}{2\alpha}}\exp(-2\alpha\rho^{\alpha}\eta/\pi + \rho x\eta)\dif\eta \\
		&\ll x^{\kappa(\Re\beta+1)+AC}\frac{\pi}{2\alpha\rho^{\alpha} - \pi\rho x}\bigl(1 - \exp(-\rho^{\alpha}+\rho x\pi/(2\alpha))\bigr)\\
		&
		 \ll x^{AC+\kappa\Re\beta-1},
	\end{aligned}
\end{multline*}
whenever $A>(\pi/(2\alpha))^{\kappa}$ so that $-\rho^{\alpha} + \rho x\pi/(2\alpha)<0$.

The remaining terms in the expression for $F_{\alpha, \beta}(z)$ can be written as
\begin{gather}
	\int_{0}^{1}t^{\beta}\biggl(\exp(\I t^{\alpha} - izt) - \sum_{\mathclap{\substack{m,n \\ m+n\alpha + \Re\beta\leq -1}}}\frac{\I^{n}(-\I z)^{m}}{n!\,m!}t^{n\alpha+m}\biggr)\dif t
	+ \sideset{}{'}\sum_{\mathclap{\substack{m,n \\ m+n\alpha + \Re\beta\leq -1}}}\frac{\I^{n}(-\I z)^{m}}{n!\,m!}\frac{1}{\beta+n\alpha+m+1} \nonumber \\
	{} + \int_{1}^{\rho}t^{\beta}\exp\bigl(\I(t^{\alpha}-xt)\bigr)\exp(yt)\dif t  \label{eq: F-L osc: rest F}.
\end{gather}

Suppose first that $\Re\beta\le-1$. By Taylor's theorem, the integrand of the first integral is bounded by $ct^{-1+\eps}\abs[1]{(-\I z)^{\lfloor-1-\Re\beta\rfloor + 1}\exp(-\I zt_{0})}$ for some constant $c$, some positive $\eps$, and some $t_{0}\in [0,1]$; hence, after integrating, this is $\ll x^{\lfloor-1-\Re\beta\rfloor + 1}$. The sum is $\ll x^{\lfloor-1-\Re\beta\rfloor}$. The last integral is $\ll x^{AC}$ if $\Re\beta<-1$ and $\ll x^{AC}\log x$ if $\Re\beta = -1$. If $\Re\beta>-1$ then we have no finite part contributions and we can integrate from 0 to $\rho$, yielding the bound $\ll_{A} x^{AC+\kappa(\Re\beta+1)}$. 

In conclusion, for $z\in \Omega_{C}$, and any fixed constant $A>\max(1, (\pi/(2\alpha))^{\kappa})$, we have
\[
	F_{\alpha, \beta}(z)\ll 
	\begin{dcases*}
		\abs{x}^{\lfloor-1-\Re\beta\rfloor+1} + \abs{x}^{AC},		 	&if $\Re\beta<-1$;\\
		 \abs{x} + \abs{x}^{AC}\log x,							&if $\Re\beta=-1$;\\
		\abs{x}^{AC+\frac{\Re\beta+1}{\alpha-1}},					&if $\Re\beta>-1$.
	\end{dcases*}
\]
\begin{remark}We end this section with two remarks.

\begin{enumerate} 	
\item [(i)]One can get better bounds on the last integral in \eqref{eq: F-L osc: rest F} by using the stationary phase principle, instead of bounding trivially. For example, when $\beta=0$ one can obtain $F_{\alpha,0}(z) \ll \abs{x}^{AC+\frac{-\alpha/2+1}{\alpha-1}}$.
	\item [(ii)]The function $\tau(x) = \exp\bigl(\I x^{1+1/\kappa}/\log^{1/\kappa}x\bigr)$ considered in the Introduction has entire Fourier transform $\hat{\tau}$, as can be shown in the same way as in Section \ref{sec: F-L osc: analytic continuation}. Similarly, one may deduce polynomial bounds for $\hat{\tau}$ in the region
	$\{ z:  \abs{\Im z} \le C(1+\abs{\Re z})^{-\kappa} \}$, by choosing $\rho \coloneqq Ax^{\kappa}\log x$ in the above procedure. 
\end{enumerate}
\end{remark}


\begin{thebibliography}{99}

\bibitem{BattyBorichevTomilov2016} C.~J.~K.~Batty, A.~Borichev, Y.~Tomilov, \emph{$L^{p}$-tauberian theorems and $L^{p}$-rates for energy decay,} J. Funct. Anal. \textbf{270} (2016), 1153--1201.


\bibitem{BattyDuyckaerts2008} C.~J.~K.~Batty, T.~Duyckaerts, \emph{Non-uniform stability for bounded semi-groups in Banach spaces,} J. Evol. Equ. \textbf{8} (2008), 765--780.

\bibitem{BorichevTomilov2010} A.~Borichev, Y.~Tomilov, \emph{Optimal polynomial decay of functions and operator semigroups,} Math. Ann. \textbf{347} (2010), 455--478. 

\bibitem{Broucke-Debruyne-VindasW-I2020} F.~Broucke, G.~Debruyne, J.~Vindas, \emph{On the absence of remainders in the Wiener-Ikehara and Ingham-Karamata theorems: a constructive approach,} Proc. Amer. Math. Soc., in press.

\bibitem{Chamizo-Ubis2007} F.~Chamizo, A.~Ubis, \emph{Some Fourier series with gaps,} J. Anal. Math. \textbf{101} (2007), 179--197.

\bibitem{ChoimetQueffelec2015} D.~Choimet, H.~Queff\'{e}lec, \emph{Twelve landmarks of twentieth-century analysis,} Cambridge University Press, New York, 2015.

\bibitem{DebruyneSchlage-PuchtaVindas2016} G.~Debruyne, J.-C.~Schlage-Puchta, J.~Vindas, \emph{Some examples in the theory of Beurling's generalized prime numbers,} Acta Arith. \textbf{176} (2016), 101--129.

\bibitem{DebruyneSeifert1} G.~Debruyne, D.~Seifert, \emph{An abstract approach to optimal decay of functions and operator semigroups,}  Israel J. Math. \textbf{233} (2019), 439--451.

\bibitem{DebruyneSeifert2} G.~Debruyne, D.~Seifert, \emph{Optimality of the quantified Ingham-Karamata theorem for operator semigroups with general resolvent growth,}  Arch. Math. (Basel) \textbf{113} (2019), 617--627.

\bibitem{DebruyneVindas2018} G.~Debruyne, J.~Vindas, \emph{Note on the absence of remainders in the Wiener-Ikehara theorem,} Proc. Amer. Math. Soc. \textbf{146} (2018), 5097--5103.

\bibitem{EstradaKanwal1990}R.~Estrada, R.~P.~Kanwal, \emph{A distributional theory for asymptotic expansions,} Proc. Roy. Soc. London Ser. A \textbf{428} (1990), 399--430. 


\bibitem {EstradaKanwalbook}R.~Estrada, R.~P.~Kanwal, \textit{A distributional approach to asymptotics. Theory and applications,} Second edition, Birkh\"{a}user, Boston, 2002.

\bibitem{ingham1935} A.~E.~Ingham, \emph{On Wiener's method in Tauberian theorems,} Proc. London Math. Soc. (2) \textbf{38} (1935), 458--480. 

\bibitem{karamata1934}J.~Karamata, \emph{\"{U}ber einen Satz von Heilbronn und Landau,} Publ. Inst. Math. (Beograd) \textbf{5} (1936), 28--38. 

\bibitem{korevaarbook} J.~Korevaar, \textit{Tauberian theory. A century of developments}, Grundlehren der Mathematischen Wissenschaften, 329, Springer-Verlag, Berlin, 2004.

\bibitem{Muger2018} M.~M\"uger, \emph{On Ikehara type Tauberian theorems with $O(x^{\gamma})$ remainders}, Abh. Math. Semin. Univ. Hambg. \textbf{88} (2018), 209--216.

\bibitem{NeytVindas2019MAE} L.~Neyt, J.~Vindas, \emph{Asymptotic boundedness and moment asymptotic expansion in ultradistibution spaces,} Appl. Anal. Discrete Math., in press, doi:10.2298/AADM191023021N.

\bibitem{PilipovicStankovicVindas2012}S.~Pilipovi\'{c}, B.~Stankovi\'{c}, J.~Vindas, \emph{Asymptotic behavior of generalized functions,} Series on Analysis, Applications and Computations, 5, World Scientific Publishing Co., Hackensack, NJ, 2012.

\bibitem{Stahn2018} R.~Stahn, \emph{Decay of $C_{0}$-semigroups and local decay of waves on even (and odd) dimensional exterior domains,} J. Evol. Equ. \textbf{18} (2018), 1633--1674. 

\bibitem{EstradaYangMAE} Y.~Yang, R.~Estrada, \emph{Asymptotic expansion of thick distributions,} Asymptot. Anal. \textbf{95} (2015), 1--19.

\end{thebibliography}
\end{document}